\documentclass[12pt,leqno]{article}
\usepackage{amssymb, amscd, amsmath, amsthm}
\usepackage{dsfont}
\allowdisplaybreaks

\def\beq{\begin{equation}}
\def\eeq{\end{equation}}

\theoremstyle{definition}
\newtheorem{definition}{Definition}
\newtheorem{example}{Example}
\newtheorem{remark}{Remark}

\theoremstyle{plain}
\newtheorem{theorem}{Theorem}

\newtheorem*{theorSper}{Sperner's Theorem}
\newtheorem*{trlemma}{Transfer Lemma}
\newtheorem*{theoremK}{Katona Theorem}
\newtheorem*{theoremKl}{Kleitman Theorem}
\newtheorem*{theoremBK}{Binomial Katona Theorem}

\newtheorem{claim}{Claim}
\newtheorem{lemma}{Lemma}

\newtheorem{corollary}{Corollary}
\newtheorem{proposition}{Proposition}

\numberwithin{proposition}{section}
\numberwithin{theorem}{section}
\numberwithin{definition}{section}
\numberwithin{claim}{section}
\numberwithin{lemma}{section}
\numberwithin{conjecture}{section}
\numberwithin{corollary}{section}
\numberwithin{equation}{section}
\numberwithin{example}{section}
\numberwithin{remark}{section}

\begin{document}

\title{Extremal set theory for the binomial norm}

\author{Peter Frankl,
R\'enyi Institute, Budapest, Hungary}

\date{}
\maketitle

\begin{abstract}
Let $[n] = \{1,2, \dots, n\}$ be the standard $n$-element set, $2^{[n]}$ its power set.
The size $|\mathcal F|$ of a family $\mathcal F \subset 2^{[n]}$ is the number of subsets $F \in \mathcal F$.
In the present paper we consider the more equitable measurement $\|\mathcal F\|$ where the contribution of each $F \in \mathcal F$ is $1\bigm/ {n\choose |F|}$.
In particular, $\|2^{[n]}\| = n + 1$.
Our main result is that for every positive integer $s$, $\|\mathcal F\| > (n + 1)s/s + 1$ implies the existence of $s$ pairwise disjoint members in $\mathcal F$. 
Moreover, this is best possible for all $n \geq s$.
The corresponding classical result of Kleitman is an easy consequence.
The relatively simple proofs seem to indicate that $\|\mathcal F\|$ is the right measurement for such problems.
\end{abstract}

\section{Introduction}
\label{sec:1}

Let $[n] = \{1,\dots, n\}$ be the standard $n$-element set and $2^{[n]}$ its power set.
Subsets of $2^{[n]}$ are called families.
For a family $\mathcal F \subset 2^{[n]}$, $|\mathcal F|$ denotes its size, that is, the number of distinct subsets in $\mathcal F$.

The first theorem in extremal set theory was proved by Sperner in 1928.
To state it let us say that $\mathcal F \subset 2^{[n]}$ is an \emph{antichain} if there are no distinct sets, $F, G \in \mathcal F$ satisfying $F \subset G$.

\begin{theorSper}[{\cite{S}}]
If $\mathcal F \subset 2^{[n]}$ is an antichain, then
\beq
\label{eq:1.1}
|\mathcal F| \leq {n\choose\left\lfloor\frac{n}{2}\right\rfloor}.
\eeq
\end{theorSper}

The most important generalisation of \eqref{eq:1.1} is the following

\smallskip
\noindent{\bf LYM Inequality} (cf.\ \cite{Y}, \cite{B}, \cite{L} and \cite{M}).
Suppose that $\mathcal F \subset 2^{[n]}$ is an antichain.
Then
\beq
\label{eq:1.2}
\sum_{F \in \mathcal F} 1\Bigm/ {n\choose |F|} \leq 1.
\eeq
To deduce \eqref{eq:1.1} from \eqref{eq:1.2} just note that ${n\choose |F|} \leq {n\choose \left\lfloor \frac{n}{2}\right\rfloor}$.
Sperner proved that equality holds in \eqref{eq:1.1} only if $\mathcal F = {[n]\choose \left\lfloor\frac{n}{2}\right\rfloor}$ or if $n$ is odd and $\mathcal F = {[n]\choose \frac{(n + 1)}{2}}$.
Similarly, equality holds in \eqref{eq:1.2} only if $\mathcal F = {[n]\choose i}$ for some $0 \leq i \leq n$.

The LYM inequality motivates the following definition.

\begin{definition}
\label{def:1.1}
For a family $\mathcal F \subset 2^{[n]}$ let us define its \emph{binomial norm}, $\|\mathcal F\|_n$ by
\beq
\label{eq:1.3}
\|\mathcal F\|_n = \sum_{F \in \mathcal F} 1\Bigm/ {n\choose |F|}.
\eeq
\end{definition}

When it causes no confusion we shall omit $n$ and write simply $\|\mathcal F\|$.

Often to compute $\|\mathcal F\|$, it is convenient to collect terms.
For $0 \leq i \leq n$ we set
\[
\mathcal F^{(i)} = \{F \in \mathcal F: \ |F| = i\} \ \text{ and } \ \varphi(i) = \bigl|\mathcal F^{(i)}\bigr| \Bigm/{n \choose i}.
\]
Obviously, $0 \leq \varphi(i) \leq 1$ and $\varphi(i)$ is the probability that (for fixed $i$) if we choose an $i$-subset of $[n]$ according to the uniform distribution on ${[n]\choose i}$, then it belongs to $\mathcal F$.
Let us note the easy fact
\beq
\label{eq:1.4}
\|\mathcal F\|_n = \sum_{0 \leq i \leq n} \varphi(i).
\eeq
Both from \eqref{eq:1.3} and \eqref{eq:1.4}
\beq
\label{eq:1.5}
0 \leq \|\mathcal F\|_n \leq n + 1 \ \text{ is obvious.}
\eeq

For some statements it is more convenient to use the \emph{relative measure} $\varrho_n(\mathcal F):= \|\mathcal F\|_n \bigm/(n + 1)$.
Note that
\beq
\label{eq:1.6}
0 \leq \varrho_n(\mathcal F) \leq 1.
\eeq

Note also the following simple identity valid for $0 \leq \ell \leq n$:
\beq
\label{eq:1.7}
\frac1{n + 1} \cdot \frac1{{n\choose \ell}} = \frac1{n + 2} \cdot \frac1{{n + 1\choose \ell}} + \frac1{n + 2} \cdot \frac1{{n + 1\choose \ell + 1}}.
\eeq
For a family $\mathcal F \subset 2^{[n]}$ one can define its natural extension in $2^{[n + 1]}$ by $\mathcal E(\mathcal F) = \bigl\{ E \subset [n + 1] : E \cap [n] \in \mathcal F\bigr\}$.

\setcounter{claim}{1}
\begin{claim}
\label{cl:1.2}
\beq
\label{eq:1.8}
\varrho_{n + 1}(\mathcal E(\mathcal F)) = \varrho_n(\mathcal F).
\eeq
\end{claim}

\begin{proof}
Simply note that for every $F \in \mathcal F$, denoting $|F|$ by $\ell$, \eqref{eq:1.7} implies
\[
\varrho_n(\{F\}) = \varrho_{n + 1} \bigl(\{F, F\cup \{n + 1\}\}\bigr).
\]
Hence \eqref{eq:1.8} holds.
\end{proof}

The central property for families that we consider in this paper is \emph{cross-dependence}.

\setcounter{definition}{2}
\begin{definition}
\label{def:1.3}
The families $\mathcal F_1, \dots, \mathcal F_s$ are called \emph{cross-dependent} if there is no choice of $F_1 \in \mathcal F_1, \dots, F_s \in \mathcal F_s$ such that $F_1, \dots, F_s$ are pairwise disjoint.
\end{definition}

\setcounter{example}{3}
\begin{example}
\label{ex:1.4}
Let $1 \leq k \leq n$ and $Y \in {[n]\choose k}$.
Let $j_1, \dots, j_s$ be integers, $0 \leq j_i \leq k$, $j_1 + \dots + j_s = k + 1$.
Define
\[
\mathcal G_i = \bigl\{G \subset [n]: |G \cap Y| \geq j_i\bigr\}.
\]
\end{example}

\setcounter{proposition}{4}
\begin{proposition}
\label{pr:1.5}
The families $\mathcal G_1, \dots, \mathcal G_s$ are cross-dependent and
\beq
\label{eq:1.9}
\varrho_n(\mathcal G_i) = \frac{k + 1 - j_i}{k + 1},
\eeq
\beq
\label{eq:1.10}
\sum_{1 \leq i \leq s} \varrho_n (\mathcal G_i) = s - 1.
\eeq
\end{proposition}

Note that \eqref{eq:1.9} is an easy consequence of \eqref{eq:1.8} while \eqref{eq:1.10} follows from \eqref{eq:1.9}.
Our main result states that the above construction is optimal.

\setcounter{theorem}{5}
\begin{theorem}
\label{th:1.6}
Suppose that $\mathcal F_1, \dots, \mathcal F_s \subset 2^{[n]}$ are cross-dependent, then
\beq
\label{eq:1.11}
\varrho_n(\mathcal F_1) + \dots + \varrho_n(\mathcal F_s) \leq s - 1.
\eeq
\end{theorem}

In terms of the binomial norm \eqref{eq:1.11} can be restated as
\beq
\label{eq:1.12}
\|\mathcal F_1\|_n + \dots + \|\mathcal F_s\|_n \leq (n + 1)(s - 1).
\eeq
In the case $\mathcal F = \mathcal F_1 = \dots = \mathcal F_s$ the condition is that $\mathcal F$ contains no $s$ pairwise disjoint members.
For $s = 2$ it is equivalent to $\mathcal F$ being \emph{intersecting}, that is, $F \cap F' \neq \emptyset$ for all $F, F'\in \mathcal F$.
Erd\H{o}s, Ko and Rado \cite{EKR} noted that not both of $F$ and $[n]\setminus F$ can be members of an intersecting family.
This implies
\beq
\label{eq:1.13}
|\mathcal F| \leq \frac12 2^n = 2^{n - 1}.
\eeq

Since ${n\choose \ell} = {n \choose n - \ell}$, the same argument implies
\beq
\label{eq:1.14}
\varrho_n(\mathcal F) \leq \frac12 .
\eeq
Theorem \ref{th:1.6} implies the following more general statement.

\setcounter{corollary}{6}
\begin{corollary}
\label{cor:1.7}
Suppose that $\mathcal F \subset 2^{[n]}$ contains no $s$ pairwise disjoint members.
Then
\beq
\label{eq:1.15}
\varrho_n(\mathcal F) \leq \frac{s - 1}{s}.
\eeq
\end{corollary}

One can see that \eqref{eq:1.15} is best possible for \emph{all} $n \geq s - 1$.
This is in great contrast to the problem of estimating $\max |\mathcal F|$ for families without $s$ pairwise disjoint members.
Even though Kleitman \cite{Kl} determined $\max|\mathcal F|$ for $n \equiv 0$ or $-1$ $(\text{\rm mod }s)$ more than 50 years ago, the general case appears to be hopelessly difficult (cf.\ discussion in later sections).
In \cite{FK1} the maximum of $|\mathcal F_1| + \dots + |\mathcal F_s|$ was determined if $\mathcal F_1, \dots, \mathcal F_k \subset 2^{[n]}$ are cross-dependent (cf. Section~\ref{sec:5}).

For the proof of Theorem \ref{th:1.6} we are going to need \emph{multi-layered} inequalities.
These are extensions of the following simple inequality proved by Kleitman.

Let us use the notation $\mathcal F_i^{(\ell)} = \mathcal F_i \cap {[n]\choose \ell}$ and $\varphi_i(\ell) = \bigl|\mathcal F_i^{(\ell)}\bigr| \bigm/{n\choose \ell}$, $0 \leq \ell \leq n$.

\paragraph{Kleitman Inequality} (\cite{Kl}).
Suppose that $\mathcal F_1 , \dots, \mathcal F_s \subset 2^{[n]}$ are cross-dependent.
Let $k_1, \dots, k_s$ be nonnegative integers satisfying $k_1 + \dots + k_s = n$.
Then
\beq
\label{eq:1.17} 
\sum_{1 \leq i \leq s} \varphi_i(k_i) \leq s - 1.
\eeq

The following is a common generalisation of \eqref{eq:1.12} and \eqref{eq:1.17}.

\paragraph{Multi-layered inequality.}
Suppose that $n \geq r \geq 0$ are integers, $\mathcal F_1, \dots, \mathcal F_s \subset 2^{[n]}$ are cross-dependent.
Let $k_1, \dots, k_s$ be nonnegative integers satisfying $k_1 + \dots + k_s + r = n$.
Then
\beq
\label{eq:1.18}
\sum_{0 \leq j \leq r} \sum_{1 \leq i \leq s} \varphi_i(k_i + j) \leq (r + 1)(s - 1).
\eeq
To deduce \eqref{eq:1.17} from \eqref{eq:1.18}, set $r = 0$.
To deduce \eqref{eq:1.12}, let $k_1 = \dots = k_s = 0$ and $r = n$.

We are going to prove Theorem \ref{th:1.6} and \eqref{eq:1.18} in the next section along with some further generalisations.

In Section~\ref{sec:3} we prove the binomial norm version of Harper's Theorem and the Katona Theorem on $t$-intersecting families.
In Section~\ref{sec:4} we explore the relation between $|\mathcal F|$ and $\|\mathcal F\|$.
In Section~\ref{sec:5} we deduce Kleitman's bound from \eqref{eq:1.15}.

\section{The proof of the main results}
\label{sec:2}

The proof is rather simple.
In its core lies a simple probabilistic argument which we present first.
Let $k_1, k_2, \dots, k_s$ and $r$ be nonnegative integers, $n = k_1 + \dots + k_s + r$.
Let $\mathcal F_i \subset 2^{[n]}$ for $1 \leq i \leq s$.

Choose a permutation $\pi = (y_1, \dots, y_n)$ of $[n]$ uniformly at random.
Define $P_1 = \bigl\{y_1, \dots, y_{k_1}\bigr\}, \dots, P_s = \bigl\{y_{k_1 + \dots + k_{s - 1} + 1}, \dots, y_{k_1 + \dots + k_s}\bigr\}$ and $R = \bigl\{y_{n - r + 1},\dots, y_n\bigr\}$.
Define the map $\varrho : R \to [r]$ by $\varrho(y_{n - r + i}) = i$, $1 \leq i \leq r$, and $\varrho(G) = \{\varrho(y) : y \in G\}$.

Let $\mathcal F_1, \dots, \mathcal F_s \subset 2^{[n]}$.
For $1 \leq i \leq r$ let us define the family $\mathcal G_i = \mathcal G_i(\pi) \subset 2^{[r]}$ by
\[
\mathcal G_i = \bigl\{G \subset [r] : \bigl(P_i \cup \varrho^{-1}(G)\bigr) \in \mathcal F_i\bigr\}.
\]

\begin{trlemma}
\begin{itemize}
\itemsep=-2pt
\item[{\rm (i)}] If $\mathcal F_1, \dots, \mathcal F_s$ are cross-dependent, then $\mathcal G_1, \dots, \mathcal G_s$ are cross-dependent as well.

\item[{\rm (ii)}] The expected value $\gamma_i(j)$ of $\bigl|\mathcal G_i^{(j)}\bigr|\bigm/{r\choose j}$ satisfies $\gamma_i(j) = \varphi_i(k_i + j)$.
\end{itemize}
\end{trlemma}

\begin{proof}
(i) Let $G_i \in \mathcal G_i$, $i = 1,\dots, s$.
Then $\bigl(P_i \cup \varrho^{-1}(G_i)\bigr) \in \mathcal F_i$.
Noting that $P_1, \dots, P_s$ are pairwise disjoint and disjoint to $R$ as well, the cross-dependence of $\mathcal F_1, \dots, \mathcal F_s$ implies the existence of $i \neq i'$ with $\varrho^{-1}(G_i) \cap \varrho^{-1}(G_{i'}) \neq \emptyset$.
This in turn implies $G_i \cap G_{i'} \neq \emptyset$ as desired.

(ii) For any fixed $0 \leq j \leq r$ and $G \in {[r]\choose j}$ the $(k_i + j)$-element set $P_i \cup\varrho^{-1}(G)$ is a uniformly random subset of $[n]$.
Thus the probability of $P_i \cup \varrho^{-1}(G)$ being in $\mathcal F_i$ is $\varphi_i(k_i + j)$.
This proves (ii).
\end{proof}

Let us turn to the proof of \eqref{eq:1.12} and \eqref{eq:1.18}.
First note that for $r = 0$ inequality \eqref{eq:1.18} is the same as \eqref{eq:1.11}.
Let us apply induction on~$n$.
For $n = 0$ the statement amounts to saying that $\emptyset$ can be a member of at most $s - 1$ of the $s$ families.

Assuming that \eqref{eq:1.12} holds for $n$ replaced by $r$ let us prove \eqref{eq:1.18} for all pairs $(r, n)$, $r \leq n$.

Let $\mathcal F_1, \dots, \mathcal F_s \subset 2^{[n]}$ be cross-dependent and define the random permutation $\pi$ and the sets $P_1, \dots, P_s$, $R$ and the families $\mathcal G_i \subset 2^{[r]}$ as above.
Applying the induction hypothesis on \eqref{eq:1.12} to $\mathcal G_1, \dots, \mathcal G_s$ in view of the Transfer Lemma we obtain
\[
\sum_{1 \leq i \leq s} \sum_{0 \leq j \leq r} \varphi_i(k_i + j) = \sum_{1 \leq i \leq s} \sum_{0 \leq j \leq r} \gamma_i(j) \leq (s - 1)(r + 1).
\]
This proves \eqref{eq:1.18}.

Next we assume that \eqref{eq:1.18} holds for $r = n$ and prove \eqref{eq:1.12} for $n + 1$.
Assume that $\mathcal F_i \subset 2^{[n]}$, $i = 1, \dots, s$, are cross-dependent.
By symmetry we may assume that $\emptyset \notin \mathcal F_s$.
Applying \eqref{eq:1.18} with $k_1 = \dots = k_{s - 1} = 0$, $k_s = 1$, $r = n$ yields
\beq
\label{eq:2.1}
\sum_{1 \leq i \leq s} \sum_{0 \leq \ell \leq n} \varphi^{(i)} (k_i + \ell) \leq (n + 1)(s - 1).
\eeq
For $1 \leq i \leq s - 1$ we have
\[
\|\mathcal F_i\| = \sum_{0 \leq \ell \leq n} \varphi^{(i)}(\ell) + \varphi^{(i)}(n + 1) \leq 1 + \sum_{0 \leq i \leq n} \varphi^{(i)}(\ell).
\]
For $i = s$ we know that $\varphi^{(s)}(0) = 0$, i.e.,
\[
\|\mathcal F_s\| = \sum_{0 \leq \ell \leq n} \varphi^{(s)}(1 + \ell).
\]
Comparing these with \eqref{eq:2.1} yields
\[
\sum_{1 \leq i \leq s} \|\mathcal F_i\| \leq (n + 1)(s - 1) + (s - 1) = (n + 2)(s - 1) \ \text{ as desired.}\eqno{\square}
\]

For an integer $0 \leq d \leq n$ and families $\mathcal F_i \subset 2^{[n]}$, $1 \leq i \leq s$, we say that $\mathcal F_1, \dots, \mathcal F_s$ are $d$-cross-dependent if there is no choice of pairwise disjoint $F_1 \in \mathcal F_1, \dots, F_s \in \mathcal F_s$ satisfying $|F_1 \sqcup \ldots \sqcup F_s| \leq d$.

\setcounter{theorem}{1}
\begin{theorem}
\label{th:2.2}
If $\mathcal F_1, \dots, \mathcal F_s \subset 2^{[n]}$ are $d$-cross-dependent, then
\[
\sum_{1 \leq i \leq s} \|\mathcal F_i\| \leq (n + 1)s - (d + 1).
\]
Note that for $d = n$ we get back the inequality \eqref{eq:1.12}.
\end{theorem}

\begin{proof}
Define $\mathcal F_{s + 1} = \bigl\{F \subset[n] : |F| \geq n - d\bigr\}$.
Then $\|\mathcal F_{s + 1}\| = d + 1$.

\setcounter{claim}{2}
\begin{claim}
\label{cl:2.3}
$\mathcal F_1, \dots, \mathcal F_{s + 1}$ are cross-dependent.
\end{claim}

\begin{proof}[Proof of the claim]
The opposite would mean the existence of pairwise disjoint sets $F_1 \in \mathcal F_1, \dots, F_{s + 1} \in \mathcal F_{s + 1}$.
Using $|F_1\sqcup \ldots \sqcup F_{s + 1}| \leq n$ and $|F_{s + 1}| \geq n - d$ we infer $|F_1\sqcup \ldots \sqcup F_s| \leq d$, contradicting the $d$-cross-dependence property of $\mathcal F_1, \dots, \mathcal F_s$.
\end{proof}

Now apply \eqref{eq:1.12} to the $s + 1$ families $\mathcal F_1, \dots, \mathcal F_{s + 1}$ to obtain $\|\mathcal F_1\| + \dots + \|\mathcal F_{s + 1}\| \leq (n + 1)s$.
Subtracting $\|\mathcal F_{s + 1}\| = d + 1$ yields Theorem \ref{th:2.2}.
\end{proof}

Let us mention that analogous results were obtained in \cite{FK1} for $|\mathcal F_1| + \dots + |\mathcal F_s|$ under the assumption of $d$-cross-dependence.

Let us conclude this section by another application of Theorem \ref{th:1.6}.

For families $\mathcal F_1, \dots, \mathcal F_s \subset 2^{[n]}$ define the family $\mathcal F_1\square\ldots \square\mathcal F_s$ by\break
$\mathcal F_1\square \ldots \square \mathcal F_s = \bigl\{F \subset [n]:$  there exist pairwise disjoint
$F_1 \in \mathcal F_1, \dots, F_s \in \mathcal F_s, \ F = F_1\sqcup \ldots \sqcup F_s\bigr\}$.

In human language $\mathcal F_1\square \ldots \square \mathcal F_s$ consists of all the subsets that can be obtained as the disjoint union of $s$ members of the $\mathcal F_i$, one from each.

\setcounter{theorem}{3}
\begin{theorem}
\label{th:2.4}
If $\mathcal F_1, \dots, \mathcal F_s \subset 2^{[n]}$ and $\mathcal F_1$ is an up-set, then
\beq
\label{eq:2.2}
\bigl\|\mathcal F_1\square \ldots \square \mathcal F_s\bigr\| \geq \|\mathcal F_1\| + \dots +\|\mathcal F_s\| - (s - 1)(n + 1).
\eeq
\end{theorem}

\begin{proof}
Set $\widetilde {\mathcal F} = \bigl\{[n] \setminus F : F \in \mathcal F_1\square\ldots \square \mathcal F_s\bigr\}$ and $\mathcal F_{s + 1} = 2^{[n]}\bigm/\widetilde{\mathcal F}$.
Note that
\beq
\label{eq:2.3}
\|\mathcal F_{s + 1}\| = (n + 1) - \|\widetilde{\mathcal F}\| = (n + 1) - \|\mathcal F_1\square \ldots \square \mathcal F_s\|.
\eeq
We claim that $\mathcal F_1, \dots, \mathcal F_s, \mathcal F_{s + 1}$ are cross-dependent.
Suppose the contrary and choose $F_1 \in \mathcal F_1, \dots, F_{s + 1} \in \mathcal F_{s + 1}$ that are pairwise disjoint.
If $F_1\sqcup \ldots \sqcup F_{s + 1} = [n]$, then $[n]\setminus F_{s + 1} = F_1\square \ldots \square F_s$ which contradicts the definition of $\mathcal F_{s + 1}$.

If the union is smaller, then we can use that $\mathcal F_1$ is an up-set and replace $F_1$ by its superset $[n]\setminus (F_2 \sqcup \ldots \sqcup F_{s + 1}]$ to obtain the same contradiction.

Now we may apply \eqref{eq:1.12} to $\mathcal F_1, \dots, \mathcal F_{s + 1}$ and infer $\|\mathcal F_1\| + \dots + \|\mathcal F_{s + 1}\| \leq s(n + 1)$.
Or equivalently,
\[
\|\mathcal F_1\| + \dots + \|\mathcal F_s\| - (s - 1)(n + 1) \leq (n + 1) - \|\mathcal F_{s + 1}\|.
\]
In view of \eqref{eq:2.3} this is equivalent to \eqref{eq:2.2}.
\end{proof}

Let us show that there are many cases when equality holds in \eqref{eq:2.2}.

\setcounter{example}{4}
\begin{example}
\label{ex:2.5}
Let $t \geq 0$ and $T \in {[n]\choose t}$.
Suppose that $t_1, \dots, t_s$ are non-negative integers satisfying $t_1 + \dots + t_s \leq t$.
Define
\[
\mathcal F_i = \bigl\{F \subset[n] : |F \cap T| \geq t_1\bigr\}.
\]
As we have shown before
\[
\varrho(\mathcal F_i) = 1 - \frac{t_i}{t + 1}.
\]
As to $\mathcal F_1 \square\ldots \square \mathcal F_s$, it is equal to
\[
\{F \subset [n] : |F \cap T|\! \geq\! t_1 \! +\! \dots\! +\! t_s\}\ \text{ implying }\ \varrho(\mathcal F_1\square... \square \mathcal F_s)\!
=\! 1\! - \frac{(t_1\! +\! \dots \! +\! t_s)}{t + 1}.
\]
Thus $\varrho(\mathcal F_1\square \ldots \square \mathcal F_s) = \varrho(\mathcal F_1) + \dots + \varrho(\mathcal F_s) - (s - 1)$ showing that equality holds in \eqref{eq:2.2}.
\end{example}

\section{The binomial Harper Theorem and Katona Theorem}
\label{sec:3}

For two sets $F$ and $G$ let $F + G$ denote their \emph{symmetric difference} (\emph{Boolean sum}), i.e., $F + G = (F \setminus G) \sqcup (G \setminus F)$.
For a family $\mathcal F \subset 2^{[n]}$ let $\partial \mathcal F$ denote its \emph{outer boundary:}
\[
\partial\mathcal F = \{G \subset [n] : G \notin \mathcal F, \, \exists F \in \mathcal F, \, |F + G| = 1\}.
\]

The usual Harper Theorem (\cite{H}, cf. \cite{FF} or \cite{F} for a simple proof) determines $\min|\partial \mathcal F|$ as a function of $|\mathcal F|$.
For the binomial norm we have a much simpler result.

\begin{theorem}
\label{th:3.1}
Suppose that $\emptyset \in \mathcal F$, $[n] \notin \mathcal F$.
Then
\beq
\label{eq:3.1}
\|\partial \mathcal F\| \geq 1
\eeq
where the inequality is strict unless $\mathcal F = \{F \subset [n] : \, |F| \leq j\}$ for some $0 \leq j < n$.
\end{theorem}

\begin{proof}
Consider an arbitrary full chain $\emptyset = G_0 \subset \dots \subset G_n = [n]$ where $|G_i| = i$.
In view of $\emptyset \in \mathcal F$ and $[n] \notin \mathcal F$ one can find a $j$, $0 \leq j < n$ such that $G_j \in \mathcal F$, $G_{j + 1} \notin \mathcal F$.
By definition, $G_{j + 1} \in \partial \mathcal F$.
Thus we proved that in all $n!$ full chains there is at least one member $G$ of $\partial \mathcal F$.
Since the same $G$ occurs in $|G|!(n - |G|)!$ full chains
\[
\sum_{G \in \partial \mathcal F} |G|! (n - |G|)! \geq n!.
\]
Dividing both sides by $n!$ we infer
\[
\|\partial \mathcal F\| = \sum_{G \in \partial\mathcal F} 1\bigm/ {n\choose |G|} \geq 1 \ \text{ as desired.}
\]
Suppose that equality holds.
By the above argument every full chain contains \emph{exactly} one member of $\partial\mathcal F$.
That is, $\partial \mathcal F$ is an antichain.
Now equality in \eqref{eq:3.1} implies that $\partial \mathcal F = {[n]\choose j + 1}$ for some $0 \leq j < n$.
This in turn easily entails $\mathcal F = \{F \subset [n]: \, |F| \leq j\}$.
\end{proof}

\setcounter{remark}{1}
\begin{remark}
\label{rem:3.2}
One can conclude the proof without using the uniqueness for maximal antichains in the LYM inequality.
Let us simply use that exactly one member of each full chain is in $\partial \mathcal F$.
Let $G$ be such a set $|G| = j + 1$.
This implies that every $j$-subset $F$, $F \subset G$ must be in $\mathcal F$.
For such an $F$ choose arbitrarily $H$, $F \subset H \subset [n]$, $|H| = j + 1$.
We claim that $H \notin \mathcal F$.
Indeed, otherwise we could continue with the chain $F \subset H \subset G \cup H \subset \ldots$ and find a new member of $\partial \mathcal F$ containing $G$, a contradiction.
Iterating we infer by connectivity of the containment graph between ${[n]\choose j + 1}$ and ${[n]\choose j}$ that $\partial \mathcal F = {[n]\choose j + 1}$.
\end{remark}

For a nonnegative integer $w < n$ we say that $\mathcal F \subset 2^{[n]}$ is $w$\emph{-union} if $|F \cup F'| \leq w$ for all $F, F'\in \mathcal F$.
Let us recall the following classical result.

\begin{theoremK}[{\cite{Ka}}]
Let $0 \leq w < n$ and suppose that $\mathcal F \subset 2^{[n]}$ is $w$-union.
Then {\rm (i)} or {\rm (ii)} hold.

{\rm (i)} $w = 2\ell$ and
\beq
\label{eq:3.2}
|\mathcal F| \leq \sum_{0 \leq i \leq \ell} {n\choose i},
\eeq

{\rm (ii)} $w = 2\ell - 1$ and
\beq
\label{eq:3.3}
|\mathcal F| \leq \sum_{0 \leq i < \ell} {n\choose i} + {n - 1\choose \ell - 1}.
\eeq
\end{theoremK}

The original proof of Katona gives the following binomial version.

\begin{theoremBK}
Let $0 \leq w < n$ and suppose that $\mathcal F \subset 2^{[n]}$ is $w$-union.
Then {\rm (i)} or {\rm (ii)} hold.

{\rm (i)} $w = 2\ell$ and
\beq
\label{eq:3.4}
\|\mathcal F\| \leq \ell + 1,
\eeq

{\rm (ii)} $w = 2\ell - 1$ and
\beq
\label{eq:3.5}
\|\mathcal F\| \leq \ell + \frac{\ell}{n}.
\eeq
\end{theoremBK}

\begin{proof}
Let $f_i$ denote the number of $i$-sets in $\mathcal F$.
The $w$-union condition implies $f_i = 0$ for $i > w$.
Katona \cite{Ka} proves the inequalities
\[
f_i + f_{w + 1 - i} \leq {n\choose i}, \ \ \ 1 \leq i < \frac{w + 1}{2}.
\]
For the case $w = 2\ell - 1$, $i = \ell$, the Erd\H{o}s--Ko--Rado Theorem \cite{EKR} yields
\[
f_\ell \leq {n - 1\choose \ell - 1} = \frac{\ell}{n} {n\choose \ell}.
\]
Noting that ${n\choose i} < {n\choose w + 1 - i}$ for $1 \leq i < \frac{w + 1}{2}$,
$\varphi(i) + \varphi(w + 1 - \ell) \leq 2$ follows.
Summing these along with $\varphi(0) \leq 1$ gives \eqref{eq:3.4} and \eqref{eq:3.5}.
\end{proof}

As one can guess from the easy proof, the inequalities \eqref{eq:3.4} and \eqref{eq:3.5} are considerably weaker than the bounds on $|\mathcal F|$.

\section{Size versus binomial norm}
\label{sec:4}

Throughout this section let $\mathcal F \subset 2^{[n]}$ be a complex (down-set), that is, $E \subset F \in \mathcal F$ implies $E \in \mathcal F$.

The relation between $|\mathcal F|$ and $\|\mathcal F\|$ is very intricate even for complexes.
However in certain cases one can deduce best possible lower bounds on $|\mathcal F|$ in terms of $\|\mathcal F\|$.

Let $\ell$ be an integer.
Define ${[n]\choose <\ell} = \{E \subset[n], |E| < \ell\}$.
This is a complex and it will remain a complex even if we add some $\ell$-element sets.
Therefore the best we can hope for in terms of a relation between $\|\mathcal F\|$ and $|\mathcal F|$ is the following.
\beq
\label{eq:4.1}
|\mathcal F| \geq \sum_{0 \leq i < \ell} {n\choose i} + (\|\mathcal F\| - \ell) {n\choose \ell}.
\eeq

\begin{definition}
\label{def:4.1}
Fix a pair $(n, \ell)$, $n > \ell \geq 0$.
If \eqref{eq:4.1} holds for all $\mathcal F$ with $\|\mathcal F\| < n + 1$, then we say that $(n, \ell)$ is \emph{perfect}.
If \eqref{eq:4.1} holds for all $\mathcal F$ with $\|\mathcal F\| < \frac{n + 1}{2}$, then we say that $(n, \ell)$ is \emph{quasi-perfect}.
\end{definition}

The main result of the present section is the following.

\setcounter{proposition}{1}
\begin{proposition}
\label{pr:4.2}
If $n \geq 3\ell$, then $(n, \ell)$ is quasi-perfect.
Moreover, unless $n = 3\ell$ and $2\leq \ell \leq 5$ it is perfect.
\end{proposition}

Recalling the fact $\|2^{[n - 1]}\|_n = \frac{n + 1}{2}$ one observes that in proving quasi-perfectness we need to consider only complexes without $(n - 1)$-element sets $(\varphi(n - 1) = 0)$.
As $\|\mathcal F\| < n + 1$ in Definition \ref{def:4.1}, we are going to assume \emph{always} that $\varphi(n) = 0$.

The proof relies on the following inequality

\setcounter{lemma}{2}
\begin{lemma}
\label{lem:4.3}
Suppose that $n, k, \ell$ are positive integers, $n - 1 > k > \ell$, $n \geq 3\ell$.
Then
\beq
\label{eq:4.2}
\sum_{\ell \leq j \leq k} {n\choose j} > (k - \ell + 1){n\choose \ell}.
\eeq
Moreover, except for $n = 3\ell$, $2 \leq \ell \leq 5$
\eqref{eq:4.2} holds for $k = n - 1$ as well.
\end{lemma}

Let us postpone the proof of \eqref{eq:4.2} and prove Proposition \ref{pr:4.2} first.

Recall the definition of $\varphi(i) = \left| \mathcal F \cap {[n]\choose i}\right|\Bigm/ {n\choose i}$ and the monotonicity $\varphi(n) \leq \varphi(n - 1) \leq \dots \leq \varphi(0)$ proved by Sperner \cite{S}.
This permits to define $\alpha(j) = \varphi(j) - \varphi(j + 1)$ for $0 \leq j \leq n - 1$, where $\alpha(j)$ are nonnegative and by definition
\beq
\label{eq:4.3}
\varphi(j) = \alpha(j) + \alpha(j + 1) + \dots + \alpha(n - 1) \ \ (\text{here we used } \varphi(n) = 0).
\eeq

In the exceptional cases $n = 3\ell$, $2 \leq \ell \leq 5$, $\alpha(n - 1) = 0$ follows from $\varphi(n - 1) = 0$, as well.
To prove all cases simultaneously, set $q = n - 2$ in the exceptional cases (i) and $q = n - 1$ for the rest.
In view of \eqref{eq:4.3} we have
\beq
\label{eq:4.4}
\varphi(j) = \alpha(j) + \alpha(j + 1) + \dots + \alpha(q).
\eeq
Define $\mathcal F^+$ $(\mathcal F^-)$ where a member $F \in \mathcal F$ is in $\mathcal F^+$ $(\mathcal F^-)$ iff $|F| \geq \ell$ $(|F| < \ell)$, respectively.
Obviously $\mathcal F = \mathcal F^+ \sqcup \mathcal F^-$ is a partition.
Multiplying both sides of \eqref{eq:4.2} by $\alpha(k)$ we obtain:
\beq
\label{eq:4.5}
\sum_{\ell \leq j \leq k} \alpha(k) {n\choose j} > \alpha(k) (k - \ell + 1) {n\choose \ell}.
\eeq

On the other hand \eqref{eq:4.4} implies
\[
\|\mathcal F^+\| = \sum_{\ell \leq j \leq q} \varphi(j) = \sum_{\ell \leq j \leq q} \sum_{j \leq k \leq q} \alpha(k) = \sum_{\ell \leq k \leq q} (k - \ell + 1) \alpha(k).
\]

Comparing the RHS with \eqref{eq:4.5} we infer
\beq
\label{eq:4.6}
{n\choose \ell} \|\mathcal F^+\|\! < \!\sum_{\ell \leq k \leq q} \alpha(k) \sum_{\ell \leq j \leq k} {n\choose j}
\! = \!\sum_{\ell \leq j \leq q} {n\choose j} \!\sum_{j \leq k \leq q} \alpha_k =\! \sum_{\ell \leq j \leq q} \varphi_j {n\choose j} = |\mathcal F^+|.
\eeq

For $\mathcal F^-$ define the nonnegative reals $\delta(i) = 1 - \varphi(i)$, $0 \leq i < \ell$ and note that
\beq
\label{eq:4.7}
\delta(0) + \dots + \delta(\ell - 1) = \ell - \|\mathcal F^-\|.
\eeq

In view of ${n\choose \ell} > {n\choose i}$ for all $0 \leq i < \ell$,
\beq
\label{eq:4.8}
|\mathcal F^-| \geq \sum_{0 \leq i < \ell} {n\choose i} - \bigl(\ell - \|\mathcal F^-\|\bigr) {n\choose \ell}.
\eeq
Using $|\mathcal F| = |\mathcal F^+| + |\mathcal F^-|$ and $\|\mathcal F^+\| = \|\mathcal F\| - \|\mathcal F^-\|$, from \eqref{eq:4.6} and \eqref{eq:4.8} we derive
\[
|\mathcal F| \geq \sum_{0 \leq i < \ell} {n\choose i} + (\|\mathcal F\| - \ell) {n\choose \ell} \ \text{ as desired.} \eqno{\square}
\]

\begin{proof}[Proof of \eqref{eq:4.2}]
Let us first note that \eqref{eq:4.2} trivially holds for both $\ell = 0$ and $\ell = 1$.

For the case $\ell \geq 2$ we need the next formula:
\beq
\label{eq:4.9}
{n\choose \ell + 1} + {n\choose \ell + 2} = {n\choose \ell} \left[ \frac{(n - \ell)(n + 1)}{(\ell + 2)(\ell + 1)}\right].
\eeq
Let us prove

\setcounter{claim}{3}
\begin{claim}
\label{cl:4.4}
For cases {\rm (i)} $\sim$ {\rm (iii)} one has
\beq
\label{eq:4.10}
\frac{(n - \ell)(n + 1)}{(\ell + 2)(\ell + 1)} > 4.
\eeq

\hspace*{2.1mm}{\rm (i)} \ $n \geq 3\ell + 2$, \ $\ell \geq 2$,

\hspace*{1mm}{\rm (ii)} \ $n \geq 3\ell + 1$, \ $\ell \geq 4$,

{\rm (iii)} \ $n \geq 3\ell$, \ \ \ \ \ \,$\ell\geq 6$.
\end{claim}

To prove \eqref{eq:4.10} note that for fixed $\ell$ the LHS is an increasing function of $n$.
Therefore it is sufficient to check the cases $n = 3\ell + a$, $a = 2, 1, 0$.

For $n = 3\ell + 2$, $n + 1 = 3(\ell + 1)$ and $\frac{n - \ell}{\ell + 2} = 2 - \frac{2}{\ell + 2} \geq \frac32$ for $\ell \geq 2$.
These prove \eqref{eq:4.10}.

For $n =  3\ell + 1$ after rearranging \eqref{eq:4.10} is equivalent to
\[
6\ell^2 + 7\ell + 2 > 4\ell^2 + 12\ell + 8 \ \text{ or }\ 2\ell^2 - 5\ell - 6 > 0,
\]
which is true for $\ell \geq 4$.

For $n = 3\ell$, \eqref{eq:4.10} is equivalent to $2\ell^2 - 10\ell - 8 > 0$ which is true for $\ell \geq 6$.
\end{proof}

\begin{proof}[The proof of Lemma \ref{lem:4.3}]
First note that ${n\choose \ell} < {n \choose j}$ for $\ell < j < n - \ell$.
Therefore \eqref{eq:4.2} is obvious for $k \leq n - \ell$.
This takes care of the cases $\ell = 0$ and $\ell = 1$.

Let us first consider the case $n \geq 3\ell + 1$, $\ell \geq 4$.
The terms ${n\choose j}$ that exceed ${n\choose \ell}$ are ${n\choose j}$ with $\ell + 1 \leq j \leq n - \ell - 1$.
This gives $n - 2\ell - 1 \geq \ell$ choices for~$j$.
Noting that ${n\choose \ell + 1}$ and ${n\choose \ell + 2}$ are the smallest, \eqref{eq:4.9} and \eqref{eq:4.10} imply
\[
\sum_{\ell + 1 \leq j \leq n - \ell - 1} {n\choose j} > 2(n - 2\ell - 1) {n\choose \ell}.
\]
Adding ${n \choose \ell}$ and ${n\choose n - \ell}$ yields
\[
\sum_{\ell \leq j \leq n - \ell} {n\choose j} > 2(n - 2\ell){n \choose \ell} > (k - \ell + 1) {n\choose \ell} \ \text{ by } \ k \leq n - 1, \ n > 3\ell.
\]
In the case $\ell = 3$, \eqref{eq:4.10} is valid for $n \geq 11$.
Therefore the only remaining case is $n = 10$, $\ell = 3$.
Then ${10\choose 4} + {10\choose 5} + {10 \choose 6} = 5.6 {10\choose 3}$ which is less than $6{10\choose 3}$.
However ${10\choose 8} + {10\choose 9} = 55 > 0.4{10\choose 3}$, proving \eqref{eq:4.2} for $k = 8$ and $9$.
Now let us consider the case $n = 3\ell$.
By Claim \ref{cl:4.4}, \eqref{eq:4.10} is true for $\ell \geq 6$.
Consequently the above proof shows that $(3\ell, \ell)$ is perfect for $\ell \geq 6$.

Let $\ell = 2$. For $n = 7$ we have
\[
{7\choose 6} + {7\choose 5} + {7\choose 4} + {7\choose 3} = 98 > 4{7\choose 2}.
\]
For $n \geq 8$ we can use Claim \ref{cl:4.4} (i) to conclude the proof.

The only remaining cases are $\ell = 3, 4, 5$.
The following equalities conclude the proof for $\ell = 3,4$.
\begin{align*}
{9\choose 4} + {9\choose 5} &= 3 \cdot {9 \choose 3},\\
{12\choose 5} + {12\choose 6} + {12\choose 7} &= 5\frac1{15} {12\choose 4}.
\end{align*}
For $\ell = 5$ using \eqref{eq:4.9} yields
\[
{15\choose 6} + {15\choose 7} + {15\choose 8} + {15\choose 9} = \frac{2\cdot 10\cdot 16}{6\cdot 7} {15\choose 5} > 7{15\choose 5}.\qedhere
\]
\end{proof}

Let us mention that for every positive $\varepsilon$ and $\ell > \ell_0(\varepsilon)$ one can show that $(n, \ell)$ is perfect for $n > (2 + \varepsilon)\ell$.

\section{Applications and outlook}
\label{sec:5}

Let us first use \eqref{eq:1.15} and Proposition \ref{pr:4.2} to deduce the following important result of Kleitman.

\begin{theoremKl}[{\cite{Kl}}]
Let $n \geq s \geq 3$.
Suppose that $\mathcal K \subset 2^{[n]}$ contains no $s$ pairwise disjoint members.
Then {\rm (i)} and {\rm (ii)} hold.

{\rm (i)} If $n = sk - 1$, then
\beq
\label{eq:5.1}
|\mathcal K| \leq \sum_{i \geq k} {n\choose i}.
\eeq

{\rm (ii)} If $n = sk$, then
\beq
\label{eq:5.2}
|\mathcal K| \leq \sum_{i > k} {n\choose i} + \frac{s - 1}{s} {n\choose k}.
\eeq
\end{theoremKl}

\begin{proof}
Note that if $\mathcal K$ contains no $s$ pairwise disjoint members, then the up-set $\mathcal K^* = \{H \subset [n]: \exists K \in \mathcal K, K \subset H\}$ has the same property as well.
Hence when proving (i) and (ii) we may assume that $\mathcal K$ itself is an up-set.
Define $\mathcal F = 2^{[n]} \setminus \mathcal K$.
Since $\mathcal K$ is an up-set, $\mathcal F$ is a complex.
In view of \eqref{eq:1.15}, $\|\mathcal K\| \leq \frac{(s - 1)(n + 1)}{s}$.
Consequently, $\|\mathcal F\| \geq \frac{n + 1}{s}$.

In the case $n = sk - 1$ we have $\|\mathcal F\| \geq k$.
Applying Proposition \ref{pr:4.2} with $\ell = k - 1$ we infer
\[
|\mathcal F| \geq \sum_{0 \leq i < k - 1} {n\choose i} + {n\choose k - 1} = \sum_{0 \leq i \leq k - 1} {n\choose i} \ \text{ implying \eqref{eq:5.1}.}
\]
In the case $n = sk$ we have $\|\mathcal F\| \geq k + \frac1{s}$.

Applying Proposition \ref{pr:4.2} with $\ell = k$ yields
\[
|\mathcal F| \geq \sum_{0 \leq i < k} {n\choose i} + \frac1{s} {n\choose k} \ \text{ implying \eqref{eq:5.2}.} \qedhere
\]
\end{proof}

What happens in the case $n \not\equiv 0$ or $-1$ $(\text{\rm mod }s)$?
As mentioned in the introduction to find the exact value of $\max |\mathcal F|$ in general appears to be very difficult.
The case of $n \equiv -2$ $(\text{\rm mod }s)$ was solved by Quinn \cite{Q} for $s = 3$ and recently by Kupavskii and the author (cf.\ \cite{FK2}, \cite{FK3}) for all $s \geq 4$.

However for the remaining congruence classes we are still in the search of the right conjectures.

Many of the constructions fail to be of the form $\{$all subsets of size less than $k\}$ plus some $k$-element sets.
This suggests that \eqref{eq:1.15} in itself is not sufficiently strong to establish their optimality.

In a sense this suggests that for cross-dependence and for families without $s$ pairwise disjoint members the binomial norm is the right setting.

As another application let us give a short proof for a recent result of Kupavskii and the author.

\begin{theorem}[{\cite{FK1}}]
\label{th:5.1}
Suppose that $n = s\ell + r$ where $0 \leq r < s$, $s \geq 3$.
Let $\mathcal F_1, \dots, \mathcal F_s \subset 2^{[n]}$ be cross-dependent.
Then
\beq
\label{eq:5.3}
\sum_{1 \leq i \leq s} |\mathcal F_i| \leq s \sum_{\ell < j \leq n} {n\choose j} + (s - r - 1) {n\choose \ell}.
\eeq
Considering the families $\mathcal H_1 = \dots = \mathcal H_{r + 1} = \{H \subset [n] : |H| > \ell\}$ and
$\mathcal H_{r + 2} = \dots = \mathcal H_s = \{H \subset [n] : |H| \geq \ell\}$ shows that \eqref{eq:5.3} is best possible.
\end{theorem}

\begin{proof}[Proof of \eqref{eq:5.3}]
Setting $\mathcal G_i = 2^{[n]}\setminus \mathcal F_i$, $1 \leq i \leq s$, \eqref{eq:5.3} is equivalent to
\beq
\label{eq:5.4}
\sum_{1 \leq i \leq s} |\mathcal G_i| \geq s \sum_{0 \leq j < \ell} {n\choose j} + (r + 1) {n\choose \ell}.
\eeq
Since $\|\mathcal G_i\| = (n + 1) - \|\mathcal F_i\|$, \eqref{eq:1.12} yields
\[
\sum \|\mathcal G_i\| \geq n + 1.
\]
Supposing that $(n, \ell)$ is perfect, using \eqref{eq:4.1} we infer
\[
\sum_{1 \leq i \leq s} |\mathcal G_i| \geq s \sum_{0 \leq j < \ell} {n\choose j} + (n + 1 - s\ell) {n\choose \ell},
\]
proving \eqref{eq:5.3}.

For $s \geq 4$ or for $s = 3$ and $r > 0$, $n \geq 3\ell + 1$ and the perfectness of the pair $(n, \ell)$ follows from Proposition \ref{pr:4.2}.
The only remaining cases are $s = 3$, $r = 0$, $n = 3\ell$.
By the same proposition $(3\ell, \ell)$ is perfect, unless $2 \leq \ell \leq 5$.
For these four cases one can use quasi-perfectness to conclude the proof.
\end{proof}


\small
\begin{thebibliography}{999}
\itemsep=-2pt
\bibitem[B]{B} B. Bollob\'as,
On generalized graphs (English, with Russian summary), {\it Acta Math. Acad. Sci. Hungar.} {\bf 16} (1965), 447--452.

\bibitem[EKR]{EKR} P. Erd\H{o}s, C. Ko, and R. Rado,
Intersection theorems for systems of finite sets,
{\it Quart. J. Math. Oxford Ser. (2)} {\bf 12} (1961), 313--320.

\bibitem[F]{F} P. Frankl,
A lower bound on the size of a complex generated by an antichain,
{\it Discrete Math.} {\bf 76} (1989).

\bibitem[FF]{FF} P. Frankl and Z. F\"uredi,
A short proof for a theorem of Harper about Hamming spheres,
{\it Discrete Math.} {\bf 34} (1981), 311--313.

\bibitem[FK1]{FK1} P. Frankl and A. Kupavskii,
Two problems on matchings in set families -- in the footsteps of Erd\H{o}s and Kleitman,
{\it J. Comb. Th. Ser. B} (2019).

\bibitem[FK2]{FK2} P. Frankl and A. Kupavskii,
Families of sets with no matching of sizes $3$ and $4$,
{\it European Journal of Combinatorics} {\bf 75} (2019), 123--135.

\bibitem[FK3]{FK3} P. Frankl and A. Kupavskii,
Families with no $s$ pairwise disjoint sets,
{\it J. London Math. Soc.} {\bf 95} (2017), 875--894.

\bibitem[H]{H} L. H. Harper,
Optimal numberings and isoperimetric problems on graphs,
{\it J. Combinatorial Theory} {\bf 1} (1966), 385--394.

\bibitem[Ka]{Ka} Gy. Katona, Intersection theorems for systems of finite sets,
{\it Acta Math. Acad. Sci. Hungar.} {\bf 15} (1964), 329--337.

\bibitem[Kl]{Kl} D. J. Kleitman, Maximal number of subsets of a finite set
no $k$ of which are pairwise disjoint, {\it J. Combinatorial Theory}
{\bf 5} (1968), 157--163.

\bibitem[L]{L} D. Lubell,
A short proof of Sperner's lemma,
{\it J. Combinatorial Theory} {\bf 1} (1966), 299.

\bibitem[M]{M} L. D. Meshalkin,
Generalization of Sperner's theorem on the number of subsets of a finite set,
{\it Theory of Probability and its Applications} {\bf 8} (1963), 203--204.

\bibitem[Q]{Q} F. C. Quinn, 
{\it Extremal properties of intersecting and overlapping families}, Ph. D. Thesis, MIT, 1987.

\bibitem[S]{S} E. Sperner,
Ein Satz \"uber Untermengen einer endlichen Menge (German),
{\it Math. Z.} {\bf 27} (1928), no. 1, 544--548.

\bibitem[Y]{Y} K. Yamamoto,
Logarithmic order of free distributive lattice,
{\it J. Math. Soc. Japan} {\bf 6} (1954), 343--353.

\end{thebibliography}
\end{document}